\documentclass{amsart}
\usepackage{hyperref}

\newcommand*{\mailto}[1]{\href{mailto:#1}{\nolinkurl{#1}}}
\newcommand{\arxiv}[1]{\href{http://arxiv.org/abs/#1}{arXiv:#1}}
\newtheorem{theorem}{Theorem}[section]

\newtheorem{lemma}[theorem]{Lemma}
\newtheorem{proposition}[theorem]{Proposition}
\newtheorem{corollary}[theorem]{Corollary}
\newtheorem{hypothesis}[theorem]{Hypothesis}

\newcommand{\R}{{\mathbb R}}

\newcommand{\C}{{\mathbb C}}

\newcommand{\be}{\begin{equation}}
\newcommand{\ee}{\end{equation}}

\newcommand{\ti}{\tilde}

\newcommand{\I}{\mathrm{i}}

\newcommand{\im}{\mathrm{Im}}
\newcommand{\re}{\mathrm{Re}}

\newcommand{\T}{\mathbb{T}}
\newcommand{\LT}{L^2(\T_\kappa;r\sigma_\kappa)}

\DeclareMathOperator{\supp}{supp}
\newcommand{\BCa}{BC_{a}}
\newcommand{\BCb}{BC_{b}}

\newcommand{\qd}{{[1]}}
\newcommand{\dom}[1]{\mathfrak{D}\left(#1\right)}
\newcommand{\Defell}{\mathfrak{D}_\ell}
\newcommand{\Deftau}{\mathfrak{D}_\tau}
\newcommand{\eps}{\varepsilon}

\newcommand{\sig}{\sigma}

\numberwithin{equation}{section}

\begin{document}

\title{Sturm--Liouville operators on time scales}

\author[J.\ Eckhardt]{Jonathan Eckhardt}
\address{Faculty of Mathematics\\ University of Vienna\\
Nordbergstrasse 15\\ 1090 Wien\\ Austria}
\email{\mailto{jonathan.eckhardt@univie.ac.at}}

\author[G.\ Teschl]{Gerald Teschl}
\address{Faculty of Mathematics\\ University of Vienna\\
Nordbergstrasse 15\\ 1090 Wien\\ Austria\\ and International
Erwin Schr\"odinger
Institute for Mathematical Physics\\ Boltzmanngasse 9\\ 1090 Wien\\ Austria}
\email{\mailto{Gerald.Teschl@univie.ac.at}}
\urladdr{\url{http://www.mat.univie.ac.at/~gerald/}}

\thanks{J. Difference Equ. Appl. {\bf 18}, 1875--1887 (2012)}
\thanks{{\it Research supported by the Austrian Science Fund (FWF) under Grant No.\ Y330}}


\begin{abstract}
We establish the connection between Sturm--Liouville equations on time scales and Sturm--Liouville equations
with measure-valued coefficients. Based on this connection we generalize several results for
Sturm--Liouville equations on time scales which have been obtained by various authors in the past.
\end{abstract}

\maketitle

\section{Introduction}
\label{sec:int}

Time scale calculus was introduced by Hilger in 1988 as a means of unifying differential and difference
calculus. Since then this approach has had an enormous impact and developed into a new field of
mathematics (see e.g.\ \cite{bope}, \cite{bope2} and the references therein). In particular, Sturm--Liouville
equations on time scales have attracted substantial interest
(see, e.g., \cite{abw, anp, anp2, dary, dary2, dary3, eh, ep, gus2, gus3, hu, hu2, hb, ko, ry, zz} and the references therein)
since it contains both continuous Sturm--Liouville equations as well as their discrete analog, Jacobi equations, as special cases.
However, efforts to unify these two cases go at least back to the seminal work of Atkinson \cite{at} or the book by Mingarelli \cite{ming}.
The approach chosen there is at first sight somewhat different and is based on Sturm--Liouville operators
with measure-valued coefficients. In this setting, derivatives have to be understood as Radon--Nikodym
derivatives and we have shown only recently \cite{timescales} that the Hilger derivative, the natural derivative on a time
scale, can in fact be viewed as a special Radon--Nikodym derivative.

The purpose of the present paper is
to show that Sturm--Liouville equations on time scales can in fact be obtained as a special case of
Sturm--Liouville equations with measure-valued coefficients. Based on this we will show how some
recent results for Sturm--Liouville operators with measure-valued coefficients can be used to extend
several results for Sturm--Liouville operators on time scales which have been obtained by various authors in the past.

\section{Sturm--Liouville differential expressions on time scales}

To set the stage we recall a few definitions and facts from time scale calculus \cite{bope}, \cite{bope2}.
Let $\T$ be a time scale, that is, a nonempty closed subset of $\R$ and let $a$, $b$ be the left and right endpoint of $\T$;
\begin{align*}
 a = \inf\,\T \quad\text{and}\quad b=\sup\,\T.
\end{align*} 
We define the forward and backward shifts on $\R$ via
\be
\sig(t) = \begin{cases}
\inf \{ s \in\T \,|\, t < s \}, & t < b,\\
b, & t \ge b,
\end{cases}
\quad
\rho(t) = \begin{cases}
\sup \{ s \in\T \,|\, t > s \}, & t > a,\\
a, & t \le a,
\end{cases}
\ee
in the usual way. Note that $\sig$ is nondecreasing right continuous and $\rho$ is
nondecreasing left continuous. 
A point $t\in\T$ is called right scattered if $\sig(t)>t$ and left scattered if $\rho(t) < t$. 
Associated with each time scale $\T$ are two Borel measure which are defined via their distribution functions
$\sig$ and $\rho$ (this procedure is standard and we refer to, e.g., \cite[Sect.~A.1]{tschroe} for a brief
and concise account). For notational simplicity we denote these measures by the same letters $\sig$ and $\rho$.
Furthermore,  we define
\begin{align*}
 \T_\kappa = \begin{cases}
               \T\backslash\lbrace a\rbrace, & \text{if }a\text{ is right scattered}, \\
               \T, & \text{else},
             \end{cases}
\quad  \T^\kappa = \begin{cases}
               \T\backslash\lbrace b\rbrace, & \text{if }b\text{ is left scattered}, \\
               \T, & \text{else}.
             \end{cases}
\end{align*}
The forward shift on the time scale $\T_\kappa$ is denoted by $\sigma_\kappa$, i.e.
\begin{align*}
 \sigma_\kappa(x) = \begin{cases}
                     \sigma(x), & \text{if }x>a, \\
                     \sigma\left(a\right), & \text{if } x\leq a.
                    \end{cases}
\end{align*}
Again for notational simplicity the corresponding measure is also denoted by $\sigma_\kappa$.

Some function $f$ on $\T$ is said to be $\Delta$-differentiable at some point $t\in\T$ if there is a 
number $f^\Delta(t)$ such that for every $\eps>0$ there is a neighborhood $U \subset\T$ of $t$ such that
\be
|f(\sig(t)) - f(s) - f^\Delta(t) (\sig(t)-s)| \le \eps |\sig(t)-s|, \quad s\in U.
\ee
Analogously one may define the notion of $\nabla$-differentiability of some function using the backward shift $\rho$. One can
show (\cite[Thm.~2.1]{gus2})
\be\label{condelnab}
f^\Delta(t) = f^\nabla(\sig(t)), \qquad f^\nabla(t) = f^\Delta(\rho(t))
\ee
for continuously differentiable functions.

A function $f$ on $\T$ is locally $\Delta$-absolutely continuous if $f$ is $\Delta$-absolutely continuous on each bounded subinterval 
 $[\alpha,\beta]_\T = [\alpha,\beta]\cap \T$ with $\alpha$, $\beta\in\T$. 
 For the notion of $\Delta$-absolute continuity on time scales, we refer to~\cite{cabviv}.
 If $f$ is locally $\Delta$-absolutely continuous then $f$ is continuous on $\T$, $f^\Delta$ exists almost everywhere on $\T^\kappa$ with respect to 
 $\sigma$ and $f^\Delta\in L^1_{loc}(\T^\kappa;\sigma)$, i.e.~$f^\Delta$ is integrable over each bounded subinterval 
 of $\T^\kappa$ with respect to $\sigma$.
Furthermore, note that if $\T$ is a bounded time scale, then $f$ is $\Delta$-absolutely continuous on $\T$, i.e.~$f^\Delta\in L^1(\T^\kappa;\sigma)$. Of course, similar statements are true for locally $\nabla$-absolutely continuous functions on time scales.

Now let $r$, $p$, $q$ be complex-valued functions on $\T_\kappa$ such that 
 $r$, $q\in L^1_{loc}(\T_\kappa;\sigma)$ and $1/p\in L^1_{loc}(\T_\kappa;\rho)$.
Consider the differential expression
\be
 \ell u = \frac{1}{r} \left(-\left(pu^\nabla\right)^\Delta + q u\right), \quad u\in\Defell,
\ee
where $\Defell$ is the maximal set of functions for which this expression makes sense.
 It consists of all functions $u$ on $\T$ which are locally $\nabla$-absolutely continuous and for which 
 $pu^\nabla$ is locally $\Delta$-absolutely continuous on $\T_\kappa$.
Consequently, we have $\ell u\in L^1_{loc}(\T_\kappa;|r|\sigma_\kappa)$ for $u\in\Defell$, where $|r|\sigma_\kappa$ is the measure $\sigma_\kappa$ weighted with the function $|r|$. 
Also note that for $u\in\Defell$, $u^\nabla$ is only locally integrable with respect to $\rho$, whereas the function
\begin{align}
 u^\qd(x) = p(x)u^\nabla(x), \quad x\in\T^\kappa,
\end{align}
is locally $\Delta$-absolutely continuous. Therefore we will mainly work with this function, referred to as the first quasi-derivative of $u$.

We remark that in particular in the case $p=r\equiv 1$ some authors also work with
\begin{align}
\ti{\ell} u = - u^{\Delta \Delta} + q\, u^\sig, \qquad u^\sig=u \circ \sig,
\end{align}
which is equivalent to the formulation we have chosen here by virtue of \eqref{condelnab}. However, this formulation does not
play nice with our scalar product and does not give rise to a self-adjoint operator in general as pointed out in \cite[Sect.~5]{dary}.

In the remaining part of this section we will show that the differential expression $\ell$ may be identified with a Sturm--Liouville differential expression on $\R$, whose coefficients are measures.
Given a function $f$ on $\T_\kappa$ we define its extension to all of $\R$ by
\begin{align}
  \bar{f}(t) = \begin{cases} f(t), & t\in\T_\kappa, \\ f(\sigma_\kappa(t)), & t\not\in\T_\kappa. \end{cases}
\end{align}
Now we are able to define locally finite complex Borel measures on $\R$ by 
\begin{align}
 \varrho(B) = \int_{B\cap\T_\kappa} r(t)d\sigma_\kappa(t), \quad \varsigma(B) = \int_B \frac{1}{\overline{p}(t)} dt \quad\text{and}\quad \chi(B) = \int_{B\cap\T_\kappa} q(t)d\sigma_\kappa(t),
\end{align}
for each Borel set $B$. Let $\Deftau$ be the set of all locally absolutely continuous functions $f$ on $\R$ such that $\overline{p}f'$ is locally absolutely continuous with respect to $\sigma_\kappa$ and consider the differential expression
\begin{align}
 \tau f = -\frac{d}{d\varrho}\frac{d f}{d\varsigma} + \frac{q}{r}f, \qquad f\in\Deftau.
\end{align}
Here the derivatives have to be interpreted as Radon--Nikodym derivatives, hence the right-hand side exists almost everywhere on $\T_\kappa$ with respect to $|r|\sigma_\kappa$. 
Consequently, we have $\tau f\in L^1_{loc}(\R;|r|\sigma_\kappa)$ for every $f\in\Deftau$.
Next we want to show that one may identify the spaces $\Defell$ and $\Deftau$ and the corresponding differential expressions.
Therefore, for each function $u\in\Defell$ let $\hat{u}$ be the extension
\begin{align}\label{defuhat}
  \hat{u}(x) = \begin{cases}
               u(x), & \text{if }x\in\T, \\
               u(\sigma_\kappa(x)) + u^\nabla(\sigma_\kappa(x)) \left(x - \sigma_\kappa(x)\right), & \text{if }x\not\in\T, \\            
               \end{cases}
\end{align}
to all of $\R$.

\begin{lemma}\label{lemTSMSL1}
 If $u\in\Defell$ then $\hat{u}\in\Deftau$ with
\begin{align}
 \overline{p}\hat{u}'(x) = \overline{pu^\nabla}(x),
\end{align}
for almost all $x\in\R$ with respect to Lebesgue measure and
\begin{align}
 \frac{d \overline{p}\hat{u}'}{d\sigma_\kappa}(x) %= \frac{d \overline{p u^\nabla}}{d\sigma|_{\T_\kappa}}(x) 
                     = \left(pu^\nabla\right)^\Delta(x)
\end{align}
for almost all $x\in\R$ with respect to $\sigma_\kappa$.
\end{lemma}

\begin{proof}
If $u\in\Defell$ then~\cite[Corollary~3.1]{cabviv} shows that $\hat{u}$ is locally absolutely continuous 
 and~\cite[Lemma~4.1]{cabviv} shows that the first equality holds.
Since by assumption $pu^\nabla$ is locally $\Delta$-absolutely continuous on $\T_\kappa$, the results 
 in~\cite[Theorem~3.3]{timescales} prove that $\overline{p}\hat{u}'$ is locally absolutely continuous with
 respect to $\sigma_\kappa$ and also that the second equality holds.
\end{proof}

The next lemma shows that extending functions of $\Defell$ yields all functions of $\Deftau$. 

\begin{lemma}\label{lemTSMSL2}
If $f\in\Deftau$ then the restriction $f|_\T$ lies in $\Defell$ with $\widehat{f|_\T}=f$.
In particular, the map \eqref{defuhat} is a bijection between $\Defell$ and $\Deftau$.
\end{lemma}

\begin{proof}
Since $f$ is continuous and linear outside of $\T$, we see that $f$ really is the extension of $f|_\T$, as claimed. Hence 
 because of~\cite[Corollary~3.1]{cabviv}, the restriction $f|_\T$ is locally $\nabla$-absolutely continuous.
 Furthermore~\cite[Lemma~4.1]{cabviv} shows that
\begin{align*}
 \overline{f|_\T^\nabla}(x) = f'(x) 
\end{align*}
for almost all $x\in\R$. Now since $\overline{pf|_\T^\nabla} = \overline{p}f'$ is locally absolutely continuous 
 with respect to $\sigma_\kappa$, the results in~\cite[Theorem~3.3]{timescales} show that 
 $pf|_\T^\nabla$ is locally $\Delta$-absolutely continuous on $\T_\kappa$ and hence $f|_\T\in\Defell$. 
\end{proof}

Now from these lemmas one sees that the time scale Sturm--Liouville differential expression $\ell$ is 
 essentially equal to the Sturm--Liouville expression $\tau$, which's coefficients are measures.

\begin{theorem}\label{thmDiffExpressions}
 For each $u\in\Defell$ and $f\in\Deftau$ we have
\begin{align*}
 \ell u(x) = \tau \hat{u}(x) \quad\text{and}\quad \tau f(x) = \ell f|_\T(x)
\end{align*}
for almost all $x\in\R$ with respect to $|r|\sigma_\kappa$.
\end{theorem}

\begin{proof}
 For $u\in\Defell$ we get from Lemma~\ref{lemTSMSL1} that
 \begin{align*}
  \frac{d}{d\varrho} \frac{d\hat{u}}{d\varsigma} = \frac{1}{r} \frac{d}{d\sigma} \overline{p}\hat{u}' = \frac{1}{r} \left(pu^\nabla\right)^\Delta,
 \end{align*}
 almost everywhere with respect to $\sigma_\kappa$, which proves the first claim.
 The second claim immediately follows from this and Lemma~\ref{lemTSMSL2}.
\end{proof}

Using the preceding theorem and the results in~\cite[Theorem~3.1]{measureSL} one readily obtains the following existence and uniqueness theorem.

\begin{theorem}
Let $g\in L^1_{loc}(\T_\kappa;|r|\sigma_\kappa)$ and $z\in\C$. Then for each $c\in\T_\kappa$, $d_1$, $d_2\in\C$ there is a unique solution $u\in\Defell$ of
\begin{align*}
 (\ell - z)u = g, \quad\text{with}\quad u(c) = d_1 \quad\text{and}\quad u^\qd(c) = d_2.
\end{align*}
\end{theorem}

Moreover,~\cite[Theorem~3.6]{measureSL} shows that the solutions from the preceding theorem depend analytically on $z\in\C$.
Given two functions $f$, $g\in\Defell$ we define the Wronski determinant $W(f,g)$ as
\begin{align*}
 W(f,g)(x) = f(x) g^\qd(x) - f^\qd(x) g(x), \quad x\in\T_\kappa.
\end{align*}
This function is locally $\Delta$-absolutely continuous on $\T_\kappa$ with derivative
\begin{align*}
 W(f,g)^\Delta(x) = \left( g(x)\ell f(x)  - f(x) \ell g(x) \right) r(x)
\end{align*}
for almost all $x\in\T_\kappa$ with respect to $\sigma_\kappa$.
In particular if $u_1$, $u_2\in\Defell$ are two solutions of $(\ell - z)u=0$ then the Wronskian $W(u_1,u_2)$ is constant with
\begin{align*}
 W(u_1,u_2) \not= 0 \quad\Leftrightarrow\quad u_1,~u_2 \text{ linearly independent}.
\end{align*}
In order to obtain self-adjoint operators we make the following additional assumptions on our coefficients.

\begin{hypothesis}\label{hypTS}\quad 
\begin{enumerate}
 \item $r$ is positive almost everywhere with respect to $\sigma_\kappa$.
 \item $q$ is real valued almost everywhere with respect to $\sigma_\kappa$.
 \item $p$ is real-valued and nonzero almost everywhere with respect to $\rho$.
 \item\label{hypTSTfour} $\T$ consists of more than four points.
\end{enumerate}
\end{hypothesis}

Note that under these assumptions, the measures $\varrho$, $\varsigma$ and $\chi$ satisfy the requirements of~\cite[Hypothesis~3.7]{measureSL}.
 The last assumption in Hypothesis~\ref{hypTS} is necessary in order to work in Hilbert spaces which are at least two-dimensional.

\section{Sturm--Liouville operators on bounded time scales}

In this section we consider the case when $\T$ is a bounded time scale. The results are special cases of the more general results in the next section.
However, since this case has attracted considerable interest in the past \cite{abw, anp, dary2, ep, gus2, ko}, we want to single out the corresponding results.
The most general results seem to be the ones in \cite{dary2}, where the case $p\in H^1$ and $q\in L^2$ is treated.

In order to obtain an operator in $\LT$ we first restrict $\Defell$ to the subspace
\begin{align*}
 \Defell^2 = \left\lbrace f\in\Defell \,\left|\, f,~\ell f\in\LT \right.\right\rbrace. 
\end{align*}
Since this subspace does not give rise to a self-adjoint operator we have to restrict it further.
Moreover, since the map \eqref{defuhat} acts as the identity on $\T$, the operators associated
with $\ell$ and $\tau$ are identical by virtue of Theorem~\ref{thmDiffExpressions}. Consequently,
the following theorem is an immediate consequence of the results in~\cite[Section~7]{measureSL}.

\begin{theorem}\label{thmTSSLboundsep}
 Let $\varphi_\alpha$, $\varphi_\beta\in[0,\pi)$ and suppose that 
 $\varphi_\alpha  \not=0$,  
 if $\sigma(a)$ is right scattered and 
 $p(b) \sin\varphi_\beta  \not= (b-\rho(b))\cos\varphi_\beta$,   
 if $\rho(b)$ is right scattered.
Then the linear operator $S$ in $\LT$ given by
\begin{align}\label{def:doms}
 \dom{S} = \left\lbrace f\in\Defell^2 \,\left| \begin{array}{l} 0=f(\sigma(a))\cos\varphi_\alpha-f^\qd(\sigma(a))\sin\varphi_\alpha\\
         0=f(b)\cos\varphi_\beta-f^\qd(b)\sin\varphi_\beta \end{array}\right. \right\rbrace
\end{align}
and $S f= \ell f$ for $f\in\dom{S}$ is well-defined and self-adjoint.
\end{theorem}

Note that here well-defined means that each function in $\dom{S}$ has a unique representative in $\Defell^2$ satisfying the boundary conditions.
In the cases excluded by the condition on the parameters $\varphi_\alpha$ and $\varphi_\beta$ this fails and $\dom{S}$ (as defined in \eqref{def:doms})
will no longer be dense. Explicitly, if $\sig(a)$ is right scattered and $\varphi_\alpha=0$, the boundary condition reads $f(\sig(a))=0$ and $\dom{S}$
lacks the corresponding one-dimensional subspace. Moreover, $f^\qd(\sig(a))= p(\sig(a))\frac{f(\sig(a))-f(a)}{\sig(a)-a}$ does not enter the boundary condition
and the value of $f(a)$ cannot be determined in terms of $f(\sig(a))$. Consequently, different values of $f(a)$ for one and the same element in $\LT$
will give rise to different values of $\ell f$ and $S$ becomes a multi-valued operator, which, however, is still self-adjoint (cf.\ \cite[Theorem~7.6]{measureSL}).
One can obtain a single-valued operator by removing $\sig(a)$ from the Hilbert space. In particular, this case is covered by starting with the time scale
$\T_\kappa$ from the outset. Similarly if $\rho(b)$ is left scattered.

The self-adjoint boundary conditions given in Theorem~\ref{thmTSSLboundsep} above are separate, i.e.\ they are given for each endpoint separately.
 As in the classical theory of Sturm--Liouville operators there are also coupled self-adjoint boundary conditions.

\begin{theorem}
 Let $\varphi\in[0,\pi)$, $R\in\R^{2\times2}$ with $\det R=1$ and suppose that 
 \begin{align}\label{eqnTSSLboundcop}
 p(b) R_{12}\not=(b-\rho(b))R_{22},
 \end{align}
 if $\sigma(a)$ and $\rho(b)$ are right scattered.
 Then the linear operator $S$ in $\LT$ given by
\begin{align*}
 \dom{S} = \left\lbrace f\in\Defell^2 \,\left| \left(\begin{matrix}f(b)\\f^\qd(b)\end{matrix}\right) 
      = e^{i\varphi}R\left(\begin{matrix}f(\sigma(a))\\f^\qd(\sigma(a))\end{matrix}\right)\right. \right\rbrace
\end{align*}
and $S f= \ell f$ for $f\in\dom{S}$ is well-defined and self-adjoint.
\end{theorem}

As before, the cases excluded by~\eqref{eqnTSSLboundcop} give rise to a multi-valued, self-adjoint linear operator in $\LT$.

Finally, we also mention that the results from Volkmer \cite{vol} for measure-valued Sturm--Liouville equations apply
to our situation. To this end we need the Lebesgue decomposition of the measure $d\sig$ which follows from
\cite[Thm.~5.2]{cabviv2}.

\begin{lemma}
The Lebesgue decomposition of $d\sig$ with respect to Lebesgue measure is given by
\be\label{ldsig}
d\sig(t) =  \chi_{\T}(t)dt + \sum_{t_n\in\T : \mu(t_n)>0} \mu(t_n) d\Theta(t-t_n),
\ee
where $d\Theta$ is the Dirac measure centered at $t=0$, $\chi_I$ is the characteristic function of a set,
and $\mu(t)=\sig(t)-t$ is the graininess of the time scale.
\end{lemma}

\begin{proof}
It suffices to show that the measures on the left and right-hand side of \eqref{ldsig} agree on every interval $I=[\alpha,\beta)$
whose endpoints are in $\T$. For such an interval we have $\sig(I)=\sig_-(\beta)-\sig_-(\alpha)=\beta-\alpha$. So let us turn to the
other side. Since $\T$ is closed, it can be written as $[a,b]$ minus a countable union of disjoint open intervals $(t_n,\sig(t_n))$.
Hence $\int_{[\alpha,\beta)} \chi_{\T}(t)dt$ gives $\beta-\alpha$ minus the intervals missing in $\T\cap[\alpha,\beta)$ and the sum over the Dirac
measures just makes up for this missing part.
\end{proof}

Now we are ready to show

\begin{theorem}
Let $S$ be a self-adjoint operator from Theorem~\ref{thmTSSLboundsep} and suppose $p>0$. Then $S$ has purely discrete spectrum and
if $\LT$ is infinite dimensional, the eigenvalues $E_0<E_1<E_2<\cdots$ have the following asymptotics
\be
\lim_{n\to\infty} \frac{E_n}{n^2} = \frac{\pi^2}{L^2}, \qquad L= \int_{\T} \sqrt{\frac{r(t)}{p(t)}} dt.
\ee
\end{theorem}

\begin{proof}
By our assumptions, $S$ satisfies the hypothesis of \cite[Thm~5.5]{vol} and the above claim holds with $L= GM(\varrho,\varsigma)$,
where $GM(\varrho,\varsigma)$ is the geometric mean of the two measures $\varrho$ and $\varsigma$ (see also the summary on p14).
As pointed out on p11 in \cite{vol} this geometric mean is given by $GM(\varrho,\varsigma)=
\sqrt{g} d\varsigma$, where $g=r\, p\, \chi_\T$ is the Radon--Nikodym
derivative of $d\varrho$ with respect to $d\varsigma$.
\end{proof}

The above theorem shows that the leading asymptotics comes from the continuous part of the time scale. If the time
scale has Lebesgue measure zero (i.e.\ $L=0$), the leading asymptotics will change and will no longer be captured
by the above theorem. See \cite{anp}, \cite{anp2} for some results in this direction.

\section{Weyl's alternative}

In this section we will allow our time scale $\T$ to be unbounded.
The crucial step in order to determine the associated self-adjoint operators is the classification of the endpoints into two cases
following the original ideas of Weyl. In the special case where $r=p\equiv 1$ and continuous $q$ this was first investigated in
\cite{zz} using the original approach via Weyl circles. Further results can be found in \cite{hu2}.

We say $\ell$ is in the limit-circle (l.c.)~case at the left endpoint $a$ if for each $z\in\C$ every solution of $(\ell - z)u=0$ lies in
$\LT$ near $a$, i.e.\ is square integrable near $a$ with respect to $r\sigma_\kappa$.
Furthermore, $\ell$ is said to be in the limit-point (l.p.)~case at $a$ if for each $z\in\C$ there is some solution of $(\ell-z)u=0$ which does not lie in $\LT$ near $a$.
Similarly one defines the limit-circle and limit-point cases for the right endpoint $b$.
For example note that each finite endpoint is in the l.c.~case.
Now~\cite[Theorem~5.2]{measureSL} yields Weyl's alternative.

\begin{theorem}
 At each endpoint, $\tau$ is either in the l.c.~or in the l.p.~case.
\end{theorem}

As in the case of bounded time scales we consider the subspace
\begin{align*}
 \Defell^2 = \left\lbrace f\in\Defell \,\left|\, f,~\ell f\in\LT \right.\right\rbrace,
\end{align*}
of $\Defell$.
Given two functions $f$, $g\in\Defell^2$ the limits 
\begin{align*}
 W(f,g)(a) = \mathop{\lim_{x\downarrow a}}_{x\in\T} W(f,g)(\sigma(x)) \quad\text{and}\quad 
   W(f,g)(b) = \mathop{\lim_{x\uparrow b}}_{x\in\T} W(f,g)(\sigma(x)),
\end{align*}
exist and are finite. According to~\cite[Lemma~5.6]{measureSL} it is possible to characterize the l.c.~and the l.p.~case in terms of the Wronskian at this endpoint.

\begin{lemma}
 $\ell$ is in the l.p.~case at $a$ if and only if
\begin{align*}
 W(f,g)(a) = 0, \quad f,~g\in\Defell^2.
\end{align*}
Furthermore, $\ell$ is in the l.c.~case at $a$ if and only if there is a $f\in\Defell^2$ such that
\begin{align*}
 W(f,f^\ast)(a) = 0 \quad\text{and}\quad W(f,g)(a)\not=0 \quad\text{for some }g\in\Defell^2.
\end{align*}
Similar results hold at the right endpoint $b$.
\end{lemma}

Concerning self-adjointness, the case when both endpoints are in the l.p.~case is the simplest, as~\cite[Theorem~6.2]{measureSL} shows.

\begin{theorem}\label{thmTSLPLP}
 If $\ell$ is in the l.p.~case at both endpoints, then the linear operator $S$ given by
 $Sf=\ell f$ for $f\in\dom{S}=\Defell^2$ is well-defined and self-adjoint.
\end{theorem}

In order to determine self-adjoint operators in the remaining cases we introduce functionals on $\Defell^2$, given by
\begin{align*}
 \BCa^1(f) = W(f,w_2)(a) \quad\text{and}\quad \BCa^2(f) = W(w_1,f)(a), \quad f\in\Defell^2,
\end{align*}
if $\ell$ is in the l.c.~case at $a$ and 
\begin{align*}
 \BCb^1(f) = W(f,w_2)(b) \quad\text{and}\quad \BCb^2(f) = W(w_1,f)(b), \quad f\in\Defell^2,
\end{align*}
if $\ell$ is in the l.c.~case at $b$. Here $w_1$ and $w_2$ are some real functions in $\Defell^2$ with 
\begin{align}\label{eqnTSboundcondw}
W(w_1,w_2)(a)=W(w_1,w_2)(b)=1.
\end{align}
Note that point evaluations in a finite endpoint are a special case of these functionals, as~\cite[Proposition~7.1]{measureSL} shows.

\begin{proposition}
 One may choose real $w_1$, $w_2\in\Defell^2$ satisfying~\eqref{eqnTSboundcondw} such that the corresponding functionals satisfy
 \begin{align*}
  \BCa^1(f) = f(\sigma(a)) \quad\text{and}\quad \BCa^2(f) = f^\qd(\sigma(a)), \quad f\in\Defell^2,
 \end{align*}
 if $\T$ is bounded from below and 
\begin{align*}
 \BCb^1(f) = f(b) \quad\text{and}\quad \BCb^2(f) = f^\qd(b), \quad f\in\Defell^2,
\end{align*}
if $\T$ is bounded from above.
\end{proposition}

Now the self-adjoint operators in the case when one endpoint is in the l.c.~case and the other is in the l.p.~case may be obtained from~\cite[Theorem~7.3]{measureSL}.

\begin{theorem}\label{thmTSLCLP}
 Suppose $\ell$ is in the l.c.~case at $a$, in the l.p.~case at $b$ and let $\varphi_\alpha\in[0,\pi)$ such that
  \begin{align*}
    \cos\varphi_\alpha w_2(\sigma(a)) + \sin\varphi_\alpha w_1(\sigma(a)) \not=0,
  \end{align*}
 if $a\in\R$ and $\sigma(a)$ is right scattered.
 Then the linear operator $S$ in $\LT$, given by
 \begin{align*}
  \dom{S} = \left\lbrace f\in\Defell^2 \,\left|\, \BCa^1(f)\cos\varphi_\alpha - \BCa^2(f)\sin\varphi_\alpha = 0 \right.\right\rbrace,
 \end{align*}
 and $Sf=\ell f$ for $f\in\dom{S}$ is well-defined and self-adjoint.
\end{theorem}

Again the excluded case gives rise to a multi-valued, self-adjoint linear operator.
If $\ell$ is in the l.p.~case at $a$ and in the l.c.~case at $b$, self-adjoint operators may be given similarly in terms of the functionals $\BCb^1$ and $\BCb^2$. We now turn to the case when both endpoints are in the
 l.c.~case. As in the case of bounded time scales, the self-adjoint operators may be divided into two classes. The case of separated boundary conditions may be obtained from~\cite[Theorem~7.6]{measureSL}.

\begin{theorem}\label{thmTSLCLCsep}
 Suppose $\ell$ is in the l.c.~case at both endpoints and let $\varphi_\alpha$, $\varphi_\beta\in[0,\pi)$ such that
  \begin{align*}
    \cos\varphi_\alpha w_2(\sigma(a)) + \sin\varphi_\alpha w_1(\sigma(a)) \not=0,
  \end{align*}
 if $a\in\R$ and $\sigma(a)$ is right scattered and
  \begin{align*}
    \cos\varphi_\beta w_2(b) + \sin\varphi_\beta w_1(b) \not=0,
  \end{align*}
 if $b\in\R$ and $\rho(b)$ is right scattered. Then the linear operator $S$ in $\LT$, given by
 \begin{align*}
  \dom{S} = \left\lbrace f\in\Defell^2 \,\left| \begin{array}{l} 0=\BCa^1(f)\cos\varphi_\alpha-\BCa^2(f)\sin\varphi_\alpha\\
         0=\BCb^1(f)\cos\varphi_\beta-\BCb^2(f)\sin\varphi_\beta \end{array}\right.\right\rbrace,
 \end{align*}
 and $Sf=\ell f$ for $f\in\dom{S}$ is well-defined and self-adjoint.
\end{theorem}

Now, the second class of self-adjoint operators in this case is defined via coupled boundary conditions.
The corresponding results may be found in~\cite[Theorem~7.6]{measureSL}.

\begin{theorem}\label{thmTSLCLCcoup}
 Suppose $\ell$ is in the l.c.~case at both endpoints and let $\varphi\in[0,\pi)$, $R\in\R^{2\times 2}$ 
  with $\det R=1$. Furthermore if $a$, $b\in\R$ and $\sigma(a)$, $\rho(b)$ are right scattered, set 
\begin{align*}
 \tilde{R} = \begin{pmatrix}
                    w_2^\qd(b) & -w_2(b) \\ -w_1^\qd(b) & w_1(b)
                  \end{pmatrix}^{-1}
     R \begin{pmatrix}
        w_2^\qd(\sigma(a)) & -w_2(\sigma(a)) \\ -w_1^\qd(\sigma(a)) & w_1(\sigma(a))
       \end{pmatrix},
\end{align*}
and assume that $\tilde{R}_{12}\not=0$.
  Then the linear operator $S$ in $\LT$, given by
 \begin{align*}
  \dom{S} = \left\lbrace f\in\Defell^2 \,\left| \left(\begin{matrix}\BCb^1(f)\\\BCb^2(f)\end{matrix}\right) 
      = e^{i\varphi}R\left(\begin{matrix}\BCa^1(f)\\\BCa^2(f)\end{matrix}\right)\right.\right\rbrace,
 \end{align*}
 and $Sf=\ell f$ for $f\in\dom{S}$ is well-defined and self-adjoint.
\end{theorem}

\section{Spectrum and resolvent}

In this section we will provide the resolvents of self-adjoint operators given in the preceding section.
 As the results in~\cite[Section~8]{measureSL} show, they turn out to be integral operators. We start with the case when both endpoints are in the l.c.~case.

\begin{theorem}
 Suppose $\ell$ is in the l.c.~case at both endpoints and $S$ is a self-adjoint operator as in Theorem~\ref{thmTSLCLCsep} or Theorem~\ref{thmTSLCLCcoup}.
  Then for each $z\in\rho(S)$ the resolvent $R_z$ is an integral operator 
 \begin{align*}
  R_z f(x) = \int_{\T_\kappa} G_z(x,y) f(y) r(y) d\sigma_\kappa(y), \quad x\in\T_\kappa,~ f\in\LT,  
 \end{align*}
 with some square integrable kernel $G_z$ on $\T_\kappa\times\T_\kappa$. If $u_1$, $u_2$ are two linearly independent solutions of $(\ell-z)u=0$,
  then there are coefficients $m^\pm_{ij}(z)\in\C$, $i$, $j\in\lbrace 1,2\rbrace$, such that the kernel is
  given by
 \begin{align*}
  G_z(x,y) = \begin{cases}
              \sum_{i,j=1}^2 m^+_{ij}(z) u_i(x) u_j(y), & \text{if }y\leq x, \\
              \sum_{i,j=1}^2 m^-_{ij}(z) u_i(x) u_j(y), & \text{if }y>x.
             \end{cases}
 \end{align*}
\end{theorem}

Since the resolvents are clearly Hilbert--Schmidt operators, as in~\cite[Corollary~8.2]{measureSL} we obtain some information about the spectrum in this case.

\begin{corollary}
 Suppose $\ell$ is in the l.c.~case at both endpoints and $S$ is a self-adjoint operator as in Theorem~\ref{thmTSLCLCsep} or Theorem~\ref{thmTSLCLCcoup}.
  Then $S$ has purely discrete spectrum, i.e.~$\sigma(S)=\sigma_d(S)$. Moreover
 \begin{align*}
  \mathop{\sum_{\lambda\in\sigma(S)}}_{\lambda\not=0} \frac{1}{\lambda^2} < \infty \quad\text{and}\quad 
    \dim\ker\left(S-\lambda\right)\leq 2, \quad \lambda\in\sigma(S).
 \end{align*}
\end{corollary}

If $S$ is a self-adjoint operator with separated boundary conditions as in Theorem~\ref{thmTSLPLP}, Theorem~\ref{thmTSLCLP} or Theorem~\ref{thmTSLCLCsep}, then the resolvent has a simpler form, as~\cite[Theorem~8.3]{measureSL} shows.

\begin{theorem}\label{thmTSressep}
 Suppose $S$ is a self-adjoint operator as in Theorem~\ref{thmTSLPLP}, Theorem~\ref{thmTSLCLP} or Theorem~\ref{thmTSLCLCsep}.
  Furthermore let $z\in\rho(S)$ and $u_a$, $u_b$ be non-trivial solutions of $(\ell-z)u=0$ such that 
 \begin{align*}
  u_a \begin{cases}
         \text{satisfies the boundary condition at }a\text{ if }\ell\text{ is l.c.~at }a, \\
         \text{lies in }\LT\text{ near }a\text{ if }\ell\text{ is l.p.~at }a,
      \end{cases}
 \end{align*}
 and
 \begin{align*}
  u_b \begin{cases}
         \text{satisfies the boundary condition at }b\text{ if }\ell\text{ is l.c.~at }b, \\
         \text{lies in }\LT\text{ near }b\text{ if }\ell\text{ is l.p.~at }b.
      \end{cases}
 \end{align*}
 Then the resolvent $R_z$ is given by
 \begin{align*}
  R_z f(x) %& = W(u_b,u_a)^{-1} \left(u_b(x) \int_{(a,x)_{\T_\kappa}} u_a f r d\sigma + u_a(x)\int_{[x,b)_{\T_\kappa}} u_b f r d\sigma\right) \\
           & = \int_{\T_\kappa} G_z(x,y) f(y)r(y) d\sigma_\kappa(y), \quad x\in\T_\kappa,~f\in\LT,
 \end{align*}
 where
 \begin{align*}
  G_z(x,y) = \begin{cases}
               W(u_b,u_a)^{-1} u_a(x) u_b(y), & \text{if }x< y, \\
               W(u_b,u_a)^{-1} u_a(y) u_b(x), & \text{if }x\geq y.
             \end{cases}
 \end{align*}
\end{theorem}

 Furthermore if $S$ is a self-adjoint operator as in Theorem~\ref{thmTSressep}, then~\cite[Corollary~8.4]{measureSL} shows that all eigenvalues of $S$ are simple.

\section{Weyl--Titchmarsh theory}

The associated eigenfunction expansion was considered in \cite{dary3}, \cite{gus2}, \cite{hb} for the case of bounded time scales and for semi-unbounded time scales in \cite{gus3}, \cite{hu}.
Here we will obtain it via classical Weyl--Titchmarsh theory thereby generalizing the presently best result from \cite{hu} where the case $r=1$, $p$ differentiable, and $q$ continuous is treated.
For further generalizations of Weyl--Titchmarsh theory to time scale systems see \cite{sbc}, \cite{shz} and the references therein.

In this section assume that our time scale $\T$ is bounded from below and let $S$ be a self-adjoint operator as in Theorem~\ref{thmTSLCLP} or Theorem~\ref{thmTSLCLCsep}. In particular there is some $\varphi_\alpha\in[0,\pi)$ such that the boundary condition at the point $\sigma(a)$ reads
\begin{align*}
 f(\sigma(a))\cos\varphi_\alpha - f^\qd(\sigma(a))\sin\varphi_\alpha = 0.
\end{align*}
Now for each $z\in\C$ consider the linearly independent solutions $\theta_z$, $\phi_z$ of $(\ell-z)u=0$ with the initial conditions
\begin{align*}
 \theta_z(\sigma(a)) = \phi_z^\qd(\sigma(a)) = \cos\varphi_\alpha \quad\text{and}\quad -\theta_z^\qd(\sigma(a)) = \phi_z(\sigma(a)) = \sin\varphi_\alpha.
\end{align*}
Note that the solution $\phi_z$ satisfies the boundary condition at $\sigma(a)$. 
Given these solutions one may define a function $m$ on the resolvent set of $S$ by requiring that the solutions
\begin{align*}
 \psi_z(x) = \theta_z(x) + m(z)\phi_z(x), \quad x\in\T,~z\in\rho(S),
\end{align*}
lie in $\LT$ near $b$ and satisfy the boundary condition at $b$ if $S$ is self-adjoint operator as in Theorem~\ref{thmTSLCLCsep}. The function $m$ is called the Weyl--Titchmarsh $m$-function of $S$, the solutions $\psi_z$, $z\in\rho(S)$ are called the Weyl solutions. Now the results in~\cite[Section~9 and Section~10]{measureSL} readily yield the following properties of $m$.

\begin{theorem}
The Weyl--Titchmarsh $m$-function is a Herglotz--Nevanlinna function. In particular, there is a unique Borel measure $\mu$ on $\R$ such that
\begin{align*}
 m(z) = \re(m(\I)) + \int_\R \frac{1}{\lambda-z} - \frac{\lambda}{1+\lambda^2}d\mu(\lambda), \quad z\in\C\backslash\R.
\end{align*}
\end{theorem}

The measure $\mu$ is called the spectral measure of $S$. Indeed the next theorem, obtained from~\cite[Section~10]{measureSL} justifies this name.

\begin{theorem}
The mapping $\mathcal{F}$ given by
\begin{align*}
 \mathcal{F} f(\lambda) = \mathop{\lim_{\beta\uparrow b}}_{\beta\in\T} \int_{(a,\beta]\cap\T_\kappa} \phi_\lambda(x)f(x) r(x)d\sigma_\kappa(x), \quad \lambda\in\R,~f\in\LT,
\end{align*}
where the limit exists in $L^2(\R;\mu)$, is unitary from $\LT$ onto $L^2(\R;\mu)$ and maps $S$ onto multiplication with the identity function in $L^2(\R;\mu)$. 
\end{theorem}

As a consequence of this theorem we may read off the spectrum of $S$ from the boundary behavior of the Weyl--Titchmarsh $m$-function in the usual way.

\begin{corollary}
 The spectrum of $S$ is given by
 \begin{align*}
  \sigma(S) = \supp(\mu) = \overline{\lbrace\lambda\in\R\,|\, 0<\limsup_{\epsilon\downarrow 0} \im(m(\lambda+\I\epsilon))\rbrace}.
 \end{align*}
\end{corollary}

\bigskip
\noindent
{\bf Acknowledgments.}
We thank Gusein Guseinov for helpful discussions and hints with respect to the literature.

\end{document}